\documentclass[11pt]{amsart}
\usepackage{amsmath}
\usepackage{amsfonts}
\usepackage{amssymb}
\newtheorem{thm}{Theorem}

\newtheorem{cor}[thm]{Corollary}

\newtheorem{lemma}[thm]{Lemma}

\newtheorem{conj}[thm]{Conjecture}

\theoremstyle{definition}
\newtheorem{remark}[thm]{Remark}

\newcommand{\bb}[1]{\mathbb{#1}}
\newcommand{\cl}[1]{\mathcal{#1}}

\begin{document}

\title[Syndetic Sets and Amenability]{Syndetic Sets and Amenability}

\author[V.~I.~Paulsen]{Vern I.~Paulsen}
\address{Department of Mathematics, University of Houston,
Houston, Texas 77204-3476, U.S.A.}
\email{vern@math.uh.edu}

\thanks{This research was supported in part by NSF grant DMS-0600191.}

\begin{abstract}
We prove that if an infinite, discrete semigroup has the property that
every right syndetic set is left syndetic, then the semigroup has a
left invariant mean. We prove that the weak*-closed convex hull of the
two-sided translates of every bounded function on an infinite discrete semigroup contains a constant function. Our proofs use the algebraic properties of the Stone-Cech compactification.
\end{abstract}

\maketitle


\section{The Main Theorems} 
Let $S$ be an infinite discrete semigroup. Given $f \in \ell^{\infty}(S)$
and $t \in S$ set $(t \cdot f)(s) = f^t(s)=f(st)$ and $r_t(f) = f^t$ which is called the {\em right translate of $f$ by $t$.} Similarly,
we set $(f \cdot t)(s) = f_t(s) =f(ts)$ and $l_t(f) = f_t$ which is called the {\em left translate of $f$ by $t$.} Moreover, $t \to t \cdot f$ and $t \to f \cdot t$ define associative, left and right actions of $S$ on $\ell^{\infty}(S).$

A classic result of Mitchell\cite{M} says that the weak*-closed convex
hull of the set of right(respectively, left) translates of $f$ contains a constant function for every $f \in \ell^{\infty}(S)$ if and
  only if $\ell^{\infty}(S)$
  has a left(respectively, right) invariant mean. A semigroup is
  called {\em left(respectively, right) amenable} if
  $\ell^{\infty}(S)$ has a left(respectively, right) invariant mean.
For a simpler proof of Mitchell's result, along with an extension to the topological case,
  see the paper of Granirer and Lau\cite{GL}.

M. M. Day proved\cite{D} that if a semigroup has a left invariant mean and a right invariant mean, then it has a mean that is simultaneously
left and right invariant, called an {\em invariant mean.}.  A
semigroup with an invariant mean is called {\em amenable.} By another
result of Day\cite{D}, any group that has a left or right invariant
mean has an invariant mean. However, Day also gave examples of
semigroups that had left(respectively, right) invariant means but no
right(respectively, left) invariant means.

In these notes we prove the following theorems.

\begin{thm} Let $S$ be an infinite discrete semigroup, let $f \in
  \ell^{\infty}(S),$ and let $\cl C(f)$ denote the weak*-closed,
  convex hull of the set of functions obtained by taking right and
  left translates of $f.$  Then $\cl C(f)$ contains a constant
  function.
\end{thm}

While there might be many ways to see this result, the proof that we will present uses syndetic sets and the algebra in the Stone-Cech compactification. This proof also leads to a new condition for amenability, involving syndetic sets.

Recall that a subset $B$ of a semigroup $S$ is called {\bf right syndetic}\cite{HS} if there is a finite set $\{t_1,...,t_k \}$ of $S$ such that $S= t_1^{-1}B \cup \cdots \cup t_k^{-1}B,$ where $t^{-1}B = \{ s \in S: ts \in B \}.$  A set $B$ is called {\bf left syndetic}\cite{HS} if there is a finite set $\{ t_1,..., t_k \}$ in $S$ such that $S = Bt_1^{-1} \cup \cdots \cup Bt_k^{-1},$ where $Bt^{-1} = \{ s \in S: st \in B \}.$ When authors speak of syndetic sets with no reference to right or left, they generally mean right syndetic.

\begin{thm} Let $S$ be an infinite, discrete semigroup. If every
  right (respectively, left)
  syndetic subset of $S$ is also left (respectively, right) syndetic and $f \in
  \ell^{\infty}(G),$ then the weak*-closed
  convex hull of the set of right (respectively, left) translates of $f$ contains a
  constant function.
\end{thm}

Combining this theorem with Mitchell's result, we have:

\begin{cor} Let $S$ be an infinite discrete semigroup. If every
  right (respectively, left)
  syndetic set in $S$ is left (respectively, right) syndetic, then
  there is a left (respectively, right) invariant mean on $\ell^{\infty}(S).$ If every right syndetic set is left syndetic and every left syndetic set is right syndetic, then $S$ is amenable.
\end{cor}

\begin{cor} Let $G$ be an infinite discrete group. If $G$ is not amenable, then $G$ contains a set that is right syndetic but not left syndetic.
\end{cor}

\begin{remark} In the case of a group, it is easily seen that a set
  $B$ is right (respectively, left) syndetic if and only if $B^{-1}$ is
  left(respectively, right) syndetic.  Thus, every right syndetic set is left syndetic if and only if $B^{-1}$ is right syndetic for every right syndetic set $B.$ It also follows for groups that every right syndetic set is left syndetic if and only if every left syndetic set is right syndetic. 
\end{remark}

\begin{remark} Although the property that every right syndetic set is left syndetic implies amenability for groups, this property does not characterize amenability. The free product of the group of order two with itself, $\bb Z_2 * \bb Z_2$ is amenable, but contains a right syndetic subset that is not left syndetic. Indeed, if we let the two generators of this free product be denoted $a, b,$ and let $B$ be the subset of all words beginning with $b,$ then $\bb Z_2 * \bb Z_2 = B \cup aB,$ but $B$ is not left syndetic. To see that $\bb Z_2 * \bb Z_2$ is amenable, note that the abelian subgroup generated by the word $ba$ is of index two.
\end{remark}

\section{Proofs of the Main Results}

Our proofs of the above theorems will use a number of facts about the Stone-Cech
compactification of $S,$ $\beta S,$ the products on $\beta S,$ and some knowledge of
ultrafilters.

First, recall that points in $\beta S$ can be identified with
ultrafilters. Given $p \in \beta S$ if we let $\cl N(p)$ denote the
set of open neighborhoods of $p,$ then the subsets of $S$ defined
by $\{ U \cap G: U \in \cl N(p) \}$ is the ultrafilter that determines
$p.$ 
 Conversely, if $p$ is an ultrafilter and $A \in p,$ then $A^{-}
\subseteq \beta S$ is a clopen set that is a neighborhood of $p.$

Recall that $\ell^{\infty}(S) = C(\beta S).$

Next there are two ways to define products on $\beta S.$  If $p,q \in
\beta S$ and $p = \lim_{\lambda} s_{\lambda},  q = \lim_{\mu}
t_{\mu},$ then we set $p \cdot q = \lim_{\lambda} [ \lim_{\mu}
s_{\lambda} t_{\mu}]$ and $p \diamond q = \lim_{\mu} [ \lim_{\lambda}
s_{\lambda} t_{\mu} ].$ These both define associative operations on
$\beta S$ that make it into a semigroup. Moreover, $p \cdot q$ is
continuous in the $p$ variable, while $p \diamond q$ is continuous in
the $q$ variable.  See \cite{HS}. Note that $s \cdot q = s \diamond q$
and $p \cdot s = p \diamond s,$ when $s \in S.$

These products on $\beta S$ induce associate actions of $\beta S$ on $\ell^{\infty}(S) = C(\beta S).$
Given $f \in C(\beta S)$ and $q \in \beta S,$  if we define $f_q, f^q$ by $f_q(p) = f(q \diamond p)$ and
$f^q(p) = f(p \cdot q),$ then $f_q,f^q \in C(\beta S).$ If we set $r_q(f) = q \cdot f = f^q,$ then $q \to q \cdot f$ defines an associative, left action of the semigroup $(\beta S, \cdot)$ on $\ell^{\infty}(S) = C(\beta S).$ While setting, $l_q(f) = f \diamond q= f_q,$ defines an associative, right action of the semigroup
$(\beta G, \diamond)$ on $\ell^{\infty}(S) = C(\beta S).$ 


\begin{lemma} Let $f \in \ell^{\infty}(G),$ let $s_{\mu} \in S,$ $q \in \beta S$ and let $s_{\mu} \to q.$
Then $s_{\mu} \cdot f_{\mu} \to q \cdot f$ and $f \diamond s_{\mu} \to f \diamond q$ in the
  weak*-topology.
\end{lemma}
\begin{proof} Since $t \cdot q = \lim_{\mu} ts_{\mu}$ in the topology
  on $\beta G,$ and $f$ is continuous, we have that $f(t \cdot q) =
  \lim_{\mu} f(ts_{\mu}).$  Thus, the net of elements $s_{\mu} \cdot f \in \ell^{\infty}(G)$
  converges pointwise and boundedly to $q \cdot f$ and hence converge in the
  weak*-topology. The other proof is similar.
\end{proof}

We are now ready to prove Theorem 1.  

\begin{proof} It will suffice to prove the theorem in the case that $f
  \ge 0.$  Since $\cl C(f)$ is compact in the weak*-topology, and the
  action of left and right translation is continuous, by Zorn's lemma,
  there will exist minimal non-empty weak*-compact, convex subsets
  of $\cl C(f),$ that are right and left translation
  invariant.
Our goal is to prove that if $\cl C_0$ is such a subset, then every
function in $\cl C_0$ is constant.  First note that if $f_1 \in \cl
C_0,$ then $\{ f_2 \in \cl C_0 : \|f_2\| \le \|f_1\| \}$ is
weak*-closed and convex and right-left translation invariant.  Hence,
by minimality is equal to $\cl C_0.$  This forces $\|f_1\|=M$ to be
constant on $\cl C_0.$  Also, since every function in $\cl C(f)$ is
positive, we have that $\sup_s f_1(s) = M,$ for every $f_1 \in \cl C_0.$

By \cite[Theorem~1.51]{HS}, $(\beta S, \cdot)$ contains a unique, non-empty minimal two-sided ideal denoted $K(\beta S).$
We claim that if $q \in K(\beta S),$ and $f_1 \in \cl C_0,$ then $f_1(q)
=M.$ Thus, $f_1$ attains its maximum at $q$.

Let us first show why this proves the theorem. Note that $\|f\| - \cl
C_0 = \{ \|f\| - f_1 : f_1 \in \cl C_0 \}$ is a minimal, weak*-closed
convex right-left translation invariant subset of $\cl C(\|f\| -f).$
Thus, by the above argument $\|f\| -f_1$ also achieves its maximum at $q$.
Thus $f_1$ attains its maximum and minimum at $q$ which forces the function $f_1$ to be constant.

We now return to the proof of the claim. Assume to the contrary that
$f_1(q) < M.$ Pick $m, f_1(q) < m < M$ and let $U = \{ p : f_1(p) < m
\}$ an open neighborhood of $q.$ Hence, $A = U \cap S \in q.$

Given $A \subseteq S,$ we let $b^{-1}A= \{ t \in S: bt \in A \}.$
Look at $B = \{ b \in S: b^{-1}A \in q \}.$ Note that $b \in B$ iff
$b^{-1} \cdot (U \cap G) \in q,$ iff $(b^{-1}U) \cap G \in q$ iff $q
\in b^{-1}U$ iff $b \cdot q \in U$ iff $f_1(b \cdot q) < m.$

By \cite[Theorem~4.39]{HS} the set $B$ is right syndetic. Hence,
there exists a finite set $\{t_1,...,t_k \}$ in $S$ such that $S = t_1^{-1}B
\cup \cdots \cup t_k^{-1}B.$ 

Now let $f_2(s) = \frac{1}{k} [ f_1(t_1s q) + \cdots +f_1(gt_ks q)].$  By the lemma, each function of $s, f_1(t_is q)$ is a
weak*-limit of left-right translates of $f_1$ and so is in $\cl C_0.$
Hence, $f_2 \in \cl C_0.$  But for $s \in S,$ there exists $i,$ so
that $s \in t_i^{-1}B.$ Thus, $b = t_is \in B.$ Hence,
$f_1(t_isq) = f_1(bq) <m.$  From this it follows that $f_2(s) \le
\frac{1}{k}[(k-1)M + m]$ for every $gs \in S.$  

Thus, $f_2 \in \cl C_0$
and $\|f_2\| < M,$ a contradiction.
This contradiction, completes the proof.
\end{proof}

We now are set to prove Theorem~2.

\begin{proof} As before it is enough to consider the case of a
  positive function. Preceeding as in the proof of the above theorem,
  we let $\cl C(f)$ denote the weak*-closed convex hull of the set of
  right translates of $f,$ and arguing as before take a minimal
  non-empty weak*-closed convex subset $\cl C_0$ that is invariant
  under right translations. For $f_1 \in \cl C_0$ and $q \in K(\beta S),$ we
  show that $f_1$ attains its maximum value at $q.$  

If we suppose that it does not attain its maximum at $q$ and define
$A$ and $B$ as
  in the previous proof, we have that $B$ is right syndetic and hence
  also left syndetic.  

Thus, there exists $\{t_1, ..., t_k \}$ such that $S = Bt_1^{-1} \cup
\cdots \cup Bt_k^{-1}.$  Now if $s \in S,$ then $b=st_i,$ for some $b \in
B$ and some $i.$  Hence, $f_1(st_iq)= f_1(bq) < m.$
Setting $f_2(s) = \frac{1}{k} [ f_1(st_1q) + \cdots +
f_1(st_kq)],$ we have that $f_2 \in \cl C_0$ and $\|f_2\| \le
\frac{1}{k}[(k-1)M + m] < M,$ a contradiction.

Thus, $f_1$ attains its maximum at $q$ and as before, $\|f\| - f_1$ also attains its maximum at $q,$ from which it follows that $f_1$ is constant.
\end{proof}


We believe that our techniques should shed some light on a problem in
the algebra of the Stone-Cech compactification.
Given an infinite discrete semigroup $S$ the semigroups $(\beta S,
\cdot)$ and $(\beta S, \diamond)$ each have minimal two-sided ideals,
$K(\beta S, \cdot)$ and $K(\beta S, \diamond).$  It is known
\cite[Theorem~13.41]{HS} that $K(\beta S, \cdot)^- \cap K(\beta S,
\diamond) \ne \emptyset$ and $K(\beta S, \cdot) \cap K(\beta S,
\diamond)^- \ne \emptyset,$ where $B^-$ denotes the closure of a set
in $\beta S.$ Many examples are known for which $K(\beta S, \cdot)
\cap K(\beta S, \diamond) \ne \emptyset.$ The question of whether or
not this intersection could ever be empty was raised in \cite[page
241]{HS}. This question was answered in \cite{B} where the first example of a semigroup for which the intersection
is empty is given. Currently, there is little known about how the
property that this intersection is either empty or non-empty relates
to other properties of semigroups.
We believe that
the following conjectures are true.

\begin{conj}  Let $S$ be an infinite, discrete semigroup. If
  $K(\beta S, \cdot) \cap K(\beta S, \diamond) \ne \emptyset,$ then $S$
  is an amenable semigroup.
\end{conj}

In particular, we believe that the following is true, which is
sufficient to prove the above conjecture.

\begin{conj} Let $S$ be an infinite, discrete semigroup. If
  $K(\beta S, \cdot) \cap K(\beta S \diamond) \ne \emptyset$ and $\cl
  C$ is a weak*-closed, convex subset of the cone of positive
  functions in $\ell^{\infty}(S)$ that is
  right translation invariant and is minimal among all such subsets,
  then every function in $\cl C$ attains its maximum at every point in
  this intersection.
\end{conj}

\paragraph{\bf Acknowledgments} The author is grateful to Tony Lau for
several valuable observations.


\end{document}